\magnification\magstephalf
\parindent 20pt
\parskip 2pt plus 1pt minus 1pt
\baselineskip 14pt

\font\sc=cmcsc10 
\def\Piit{{\mit\Pi}}
\def\Sigit{{\mit\Sigma}}
\def\Gamit{{\mit\Gamma}}
\def\display#1:#2:#3\par{\par\hangindent #1 \noindent
			\hbox to #1{\hfill #2 \hskip .1em}\ignorespaces#3 \par}
\def\disleft#1:#2:#3\par{\par\hangindent#1\noindent
			 \hbox to #1{#2 \hfill \hskip .1em}\ignorespaces#3\par}
\def\blackslug{\hbox{\hskip 1pt \vrule width 4pt height 6pt depth 1.5pt \hskip 1pt}}

\centerline{\bf Efficient Representation of Perm Groups}
\centerline{by Donald E. Knuth\footnote{}{%
This research
was supported in part by the National Science Foundation under grant
CCR-86-10181,
and by Office of Naval Research contract
N00014-87-K-0502.}}
\centerline{\sl Computer Science Department, Stanford University}
\centerline{Dedicated to the memory of Marshall Hall}

\medskip
{\narrower\smallskip\noindent
{\bf Abstract:} This note presents an elementary version of Sims's 
algorithm for computing strong generators of a given perm group, together
with a proof of correctness and some notes about appropriate low-level data
structures. Upper and lower bounds on the running time are also obtained.
(Following a suggestion of Vaughan Pratt, we adopt the convention that
perm $=$ permutation, perhaps thereby saving millions of syllables in
future research.)
\smallskip}

\medskip
\noindent
{\bf 1.\enspace A data structure for perm groups.}\enspace
A ``perm,'' for the purposes of this paper, is a one-to-one mapping of a set
onto itself. If $\alpha$ and $\beta$ are perms such that $\alpha$ takes
$i\mapsto j$ and $\beta$ takes $j\mapsto k$, the product $\alpha\beta$
takes $i\mapsto k$. We write $\alpha^-$ for the inverse of the perm~$\alpha$;
hence $\alpha\beta=\gamma$ iff $\alpha=\gamma\beta^-$.

Let $\Piit(k)$ be the set of all perms of the positive integers that fix all
points $>k$. Consider the following data structure:
For $1\le j\le k$, either $\sigma_{kj}=\emptyset$ or $\sigma_{kj}$ is a perm
of $\Piit(k)$ that takes $k\mapsto j$. Let $\Sigit(k)$ be the set of all
non-$\emptyset$ perms~$\sigma_{kj}$. We assume that $\sigma_{kk}$ is the
identity perm; hence $\Sigit(k)$ is always nonempty.

We write $\Gamit(k)$ for the set of all perms that can be written as products
of the form $\sigma_1\,\ldots\,\sigma_k$ where each~$\sigma_i$ is in
$\Sigit(i)$. There is an
easy way to test if a given perm $\pi\in\Piit(k)$ is a member of $\Gamit(k)$:
Let $\pi$ take $k\mapsto j$. Then if $\sigma_{kj}=\emptyset$
we have $\pi\not\in\Gamit(k)$; otherwise if $k=1$ we have
$\pi\in\Gamit(k)$; otherwise $\pi\in\Gamit(k)$ iff
$\pi\sigma_{kj}^-\in\Gamit(k-1)$.

The data structure also includes a set $T(k)\subseteq\Piit(k)$
with the invariant property that each element of $\Gamit(k)$ can be written
as a product of elements of~$T(k)$. 
In other words,
$\Gamit(k)$ will be a subset of the group $\langle T(k)\rangle$ generated
by~$T(k)$, for all~$k$, throughout the course of the algorithm to be described.
(Since all elements $\pi$ of~$\Piit(k)$ are finite perms, we have
$\pi^-=\pi^r$ for some $r>0$; hence closure under multiplication
implies closure under inversion.)

The data structure is said to be up-to-date of order~$n$ if 
$\Gamit(k)\supseteq T(k)$ and if $\Gamit(k)$ is
closed under multiplication, i.e., if $\Gamit(k)=\langle T(k)\rangle$,
for $1\le k\le n$. In that case we say that the perms $\bigcup_{k=1}^n\Sigit(k)$
form a {\sl transversal system\/} of~$\Gamit(n)$, and that
the perms $\bigcup_{k=1}^n T(k)$ are {\sl strong generators\/} of~$\Gamit(n)$.
Having a transversal system
makes it easy to determine what perms are generated by a given set
of perms~$T(n)$.

\medskip\noindent
{\bf 2.\enspace Maintaining the data structure.}\enspace
Let us now discuss two algorithms that can be used to transform the data
structure when a new perm is introduced into~$T(k)$. We will first look at
the algorithms, then discuss why they are valid.

\vfill\eject

\noindent
{\sc Algorithm} $A_k(\pi)$.\enspace
Assuming that the data structure is up-to-date of order~$k$, and that
$\pi\in\Piit(k)$ but $\pi\not\in\Gamit(k)$, this procedure appends~$\pi$
to~$T(k)$ and brings the data structure back up-to-date so that $\Gamit(k)$
will equal the new $\langle T(k)\rangle$.

\display 50pt:Step A1.:
Insert $\pi$ into the set $T(k)$.

\display 50pt:Step A2.:
Perform algorithm $B_k(\sigma\tau)$ for all $\sigma\in\Sigit(k)$ and 
$\tau\in T(k)$ such that $\sigma\tau$ is not already known to be a member
of~$\Gamit(k)$. 
(Algorithm~$B_k$ may increase the size of~$\Sigit(k)$; any
new perms~$\sigma$ that are added to~$\Sigit(k)$ must also be included
in this step. Implementation details are discussed in Section~3 
below.)\quad\blackslug

\medskip
\noindent
{\sc Algorithm} $B_k(\pi)$.\enspace
Assuming that the data structure is up-to-date of order $k-1$, and that
$\pi\in\langle T(k)\rangle$, this procedure ensures that $\pi$ is
in~$\Gamit(k)$ and that the data structure remains up-to-date of
order $k-1$. (The value of $k$ will always be greater than~1.)

\display 50pt:Step B1.:
Let $\pi$ take $k\mapsto j$.

\display 50pt:Step B2.:
If $\sigma_{kj}=\emptyset$, set $\sigma_{kj}\gets\pi$ and terminate the algorithm.

\display 50pt:Step B3.:
If $\pi\sigma_{kj}^-\in\Gamit(k-1)$, terminate the algorithm.
(This test for membership in $\Gamit(k-1)$ has been described in Section~1 above.)

\display 50pt:Step B4.:
Perform algorithm $A_{k-1}(\pi\sigma_{kj}^-)$.\quad\blackslug

\medskip
The correctness of these mutually recursive procedures follows readily from
the stated invariant relations, except for one nontrivial fact: We must
verify that $\Gamit(k)$ is closed under multiplication at the conclusion
of algorithm $A_k(\pi)$. This is obvious when $k=1$, so we may assume that
$k>1$. Let $\alpha$ and~$\beta$ be elements of $\Gamit(k)$. By definition
of $\Gamit(k)$ we can write $\alpha=\gamma\sigma$, where $\gamma\in\Gamit(k-1)$
and $\sigma\in\Sigit(k)$; and by the invariant relation $\Gamit(k)\subseteq
\langle T(k)\rangle$ we can write $\beta=\tau_1\,\ldots\,\tau_r$ where each
$\tau_i\in T(k)$. We know that $\sigma\tau_1\in\Gamit(k)$, by step~A2;
hence $\sigma\tau_1=\gamma_1\sigma_1$ for some $\gamma_1\in\Gamit(k-1)$
and some $\sigma_1\in\Sigit(k)$. Similarly $\sigma_1\tau_2=\gamma_2\sigma_2$,
etc., and we finally obtain $\alpha\beta=\gamma\,\gamma_1\,\ldots\,\gamma_r
\sigma_r$. This proves that $\alpha\beta\in\Gamit(k)$, since
$\gamma\,\gamma_1\,\ldots\,\gamma_r$ is in $\Gamit(k-1)$ by induction.

\medskip
\noindent
{\bf 3.\enspace Low-level implementation hints.}\enspace
Let $s(k)$ be the cardinality of $\Sigit(k)$ and $t(k)$ the cardinality
of~$T(k)$. The algorithms of Section~2 can perhaps be implemented most
efficiently in practice by keeping a linear list of the
perms $\tau(k,1)\,\ldots\,\tau\bigl(k,t(k)\bigr)$ of~$T(k)$,
for each~$k$, together with an array of pointers to the
representations
of each~$\sigma_{kj}$ for $1\leq j<k$, using a null pointer to
represent the relation $\sigma_{kj}=\emptyset$. It is also convenient
to have a linear list $j(k,1)\,\ldots\,j\bigl(k,s(k)\bigr)$ of the
indices of the non-$\emptyset$ perms~$\sigma_{kj}$, where $j(k,1)=k$.
We will see below that the algorithm often completes its task without
needing to make many of the sets~$\Sigma(k)$ very large; thus most of
the $\sigma_{kj}$ are often~$\emptyset$.
Pointers can be used to avoid duplications between $T(k)$ and $T(k-1)$.

There are two fairly simple ways to handle the loop over~$\sigma$
and~$\tau$ in step~A2; one is recursive and the other is iterative.
The recursive method replaces step~A2 by the following operation:
``Perform algorithm $B_k(\sigma\pi)$ for all~$\sigma$ in the current
set~$\Sigma(k)$.'' Then step~B2 is also changed: ``If
$\sigma_{kj}=\emptyset$, set $\sigma_{kj}\gets \pi$ and perform
$B_k(\pi\tau)$ for all~$\tau$ in the current set $T(k)$, then
terminate the algorithm.''

The iterative method maintains an additional table, in order to
remember which pairs $(\sigma,\tau)$ have already been tested in
step~A2. This table consists of counts $c(k,i)$ for each~$k$ and for
$1\leq i\leq s(k)$, such that the product $\sigma_{kj(k,i)}\tau(k,l)$
is known to be in $\Gamit(k)$ for $1\leq l\leq c(k,i)$. When step~B2
increases the value of~$s(k)$, the newly created count
$c\bigl(k,s(k)\bigr)$ is set to zero. Step~A2 is a loop of the form

\smallskip
\halign{\qquad\qquad#\hfil\cr
$i\gets 1$;\cr
{\bf while} $i\leq s(k)$ {\bf do}\cr
\qquad{\bf begin while} $c(k,i)<t(k)$ {\bf do}\cr
\qquad\qquad {\bf begin} $l\gets c(k,i)+1$;\cr
\qquad\qquad $B_k\bigl(\sigma_{kj(k,i)}\tau(k,l)\bigr)$;\cr
\qquad\qquad $c(k,i)\gets l$;\cr
\qquad\qquad {\bf end};\cr
\qquad $i\gets i+1$;\cr
\qquad {\bf end};\cr}

\smallskip
\noindent
the invocation of $B_k$ may increase $s(k)$, but it can change $t(k')$
and $c(k',i')$ only for values of~$k'$ that are less than~$k$.

The iterative method carries out its tests in a different order from
the recursive method, so it might yield a different traversal system.

It is convenient to represent each perm $\sigma$ of $\Sigit(k)$ indirectly
in an array~$q$ that gives inverse images, so that $\sigma$ takes
$q[i]\mapsto i$ for $1\le i\le k$. All other perms~$\pi$ can be represented
directly in an array~$p$, with $\pi$ taking $i\mapsto p[i]$ for
$1\le i\le k$. To compute the direct representation~$d$ of the 
product~$\pi\sigma^-$, we can then simply set $d[i]\gets q\bigl[p[i]\bigr]$ for
$1\le i\le k$. To compute the direct representation~$d$ of the product~$\sigma\pi$,
we set $d\bigl[q[i]\bigr]\gets p[i]$ for $1\le i\le k$. Thus, the elementary operations
are fast.

\medskip
\noindent
{\bf 4.\enspace Upper bounds on the running time.}\enspace
The ``inner loop'' of the updating algorithms occurs in step~B3, the membership
test. Testing for membership of $\pi\in\Gamit(k)$ involves multiplication
by some sequence of non-identity perms
$\sigma_{k_1j_1}^-,\ldots,\sigma_{k_rj_r}^-$,
where
$k\geq k_1>\cdots >k_r>0$;
so the running time is essentially proportional to $k+k_1+\cdots +k_r$,
which is $O(k^2)$ in the worst case.

The total number of executions of $B_k(\sigma\tau)$
is $s(k)t(k)$, and we have $s(k)\le k$. The value of $t(k)$ increases by~1
each time we perform $A_k(\pi)$; every time we do this, we increase
$\Gamit(k)$ to a larger subgroup of~$\Piit(k)$, hence $t(k)$ cannot
exceed the length of the longest chain of subgroups of the symmetric
group~$\Piit(k)$. A~straightforward upper bound is therefore
$t(k)\le \theta(k!)=O(k\log\log k)$, where $\theta(N)$ is the number of prime
divisors of~$N$ counting multiplicity. Babai~[1] has shown that
$\Piit(k)$ admits no subgroup chains of length exceeding $2k-3$, when $k\ge 2$;
hence we have the sharper estimate $t(k)=O(k)$. 

It follows that algorithm $B_k(\sigma\tau)$ is performed $O(k^2)$ times,
and each occurrence of step~B3 takes $O(k^2)$ units of time. Summing for
$1\le k\le n$ allows us to conclude that {\sl a transversal system for a perm group
generated by $m$ perms of $\Piit(n)$ can be found in at most $O(n^5)+O(mn^2)$
steps}. (The term $O(mn^2)$ comes from $m$~membership tests, which are carried
out on each generator~$\pi$ before algorithm $A_n(\pi)$ is applied.)

The storage requirement for each non-identity
perm of $\Sigit(k)$ or $T(k)$ is $O(k)$;
hence we need at most $O(k^2)$ memory cells for perms of~$\Piit(k)$, and
$O(n^3)$ memory cells in all.

\medskip
\noindent
{\bf 5.\enspace A sparse example.}\enspace
Actual computations with these procedures rarely
take as much time as our worst-case estimates predict. We can learn more about
the true efficiency by studying particular cases in detail. Let us
therefore consider first the case of a group generated by a single
non-identity perm $\pi\in\Piit(n)$.

We begin, of course, with $\sigma_{kj}=\emptyset$ for $1\le j<k\le n$ and
$T(k)=\emptyset$ for $1\le k\le n$; the data structure is then up-to-date
of order~$n$, and we can perform $A_n(\pi)$. Suppose $\pi$ takes
$n\mapsto a_1\mapsto\cdots\mapsto a_{r-1}\mapsto n$. Then $A_n(\pi)$ will set
$T(n)\gets\{\pi\}$ and $\sigma_{na_j}\gets\pi^j$ for $1\le j<r$, and it
will invoke $A_{n-1}(\pi^r)$ (unless $\pi^r$ is the identity perm, in which
case the algorithm will terminate).

If, for example, we have
$$\pi=[1,2,3,4,5,6,7,14]\,[8,9,10,13]\,[11,12]$$
in cycle form, the algorithm will set $\sigma_{14,j}\gets\pi^j$ for
$1\le j<8$, and it will terminate with $T(14)=\{\pi\}$ and with all
other $T(k)$ empty. But if we relabel points 12 and~14, obtaining the
conjugate perm
$$\bar\pi=[1,2,3,4,5,6,7,12]\,[8,9,10,13]\,[11,14]\,,$$
the algorithm will act quite differently: The nontrivial perms
$\sigma_{kj}$ and sets $T(k)$ will now be
$$\eqalign{\sigma_{14,11}&=\bar\pi\vphantom{^1}\,,\cr
   T(14)&=\{\bar\pi\vphantom{^1}\}\,,\cr}\qquad
\eqalign{\sigma_{13,9}&=\bar\pi^2\,,\cr T(13)&=\{\bar\pi^2\}\,,\cr}\qquad
\eqalign{\sigma_{12,4}&=\bar\pi^4\,;\cr T(12)&=\{\bar\pi^4\}\,.\cr}$$

When the algorithm terminates, it has produced a transversal system
by which we can test if a given perm $\rho$ is a power of $\pi$ or $\bar\pi$,
respectively. In the first case this membership test involves at most one
multiplication, by $\sigma_{14,j}$ if $\rho$ takes $14\mapsto j$ where $j<8$.
In the second case the test will involve three multiplications if we have,
say, $\rho=\bar\pi^7$.

These perms $\pi$ and $\bar\pi$ are the special case $h=4$ of an infinite
family of perms of degree $n=2^h-2$, having cycles of lengths $2^{h-1}$,
$2^{h-2}$, \dots,~$2^1$. In general $\pi$ will cause $\sim{1\over2}n$
slots $\sigma_{kj}$ to become nonempty, and it will terminate after performing
$\sim{1\over2}n^2$ elementary machine steps, yielding a membership
test whose worst-case running time is $\sim n$. The corresponding perm
$\bar\pi$ will cause only $\sim\lg n$ slots $\sigma_{kj}$ to become
nonempty, and it will terminate after $\sim 2n\lg n$ steps,
yielding a membership test whose worst-case running
time is $\sim n\lg n$. Thus, the algorithm's performance can change
substantially when only two points of its input perm are relabeled.

\medskip
\noindent
{\bf 6.\enspace A dense example.}\enspace
 The algorithm needs to work harder when we wish
to find the group generated by $\{\pi_2,\pi_3,\ldots,\pi_n\}$, where
$\pi_k\in\Piit(k)$ takes $k\mapsto k-1$, and where the generators $\pi_k$
are input in increasing order of $k$. Then it is not difficult to verify
by induction that the algorithm will terminate with 
$T(k)=\{\pi_2,\ldots,\pi_k\}$
and with $\sigma_{kj}\neq\emptyset$ for $1\leq j<k\leq n$.
 Thus, the algorithm will fill all of the slots $\sigma_{kj}$,
thereby implicitly deducing that each $\Gamit(k)$ is the full
symmetric group~$\Piit(k)$.

Moreover, if the recursive method of Section 3 is being used to
implement step~A2, the algorithm will terminate with
$$\sigma_{kj}=\pi_k\pi_{k-1}\ldots\pi_{j+1}\,,
\qquad\hbox{for $1\le j<k\le n$}.$$
For after $\sigma_{kj}$ is defined, the modified step B2 will continue
to test whether the perms
$$\sigma_{kj}\pi_2\,,\quad \sigma_{kj}\pi_3\,,\quad\ldots\,,\quad
\sigma_{kj}\pi_k$$
belong to the current $\Gamit(k)$. The first $j-2$ tests will succeed;
then $B_k(\sigma_{kj}\pi_j)$ will cause $\sigma_{k,j-1}$ to be
defined. And by the time the recursive call on $B_k(\sigma_{kj}\pi_j)$
returns control to $B_k(\sigma_{kj})$, the values of~$\sigma_{ki}$
will be non-$\emptyset$ for all $i<k$; hence the remaining tests on
$\sigma_{kj}\pi_l$ for $l>j$ will succeed.

Let us examine the special case of this construction in which each $\pi_k$
is the simple transposition $[k,\,k-1]$. How much time is taken by the
$\Theta(n^3)$ membership tests $\sigma_{kj}\pi_i\in\Gamit(k)$? We have
$$\sigma_{kj}=[j,\,j+1,\,\ldots\,,\,k]\,,$$
and it follows that
$$\sigma_{kj}\pi_i=\cases{\sigma_{i,i-1}\,\sigma_{kj}\,,&if $1<i<j$;\cr
 \sigma_{ki}\,,&if $i=j+1$;\cr
 \sigma_{i-1,i-2}\,\sigma_{kj}\,,&if $i>j+1$.\cr}$$
Each membership test therefore involves at most two multiplications by
non-identity perms, and the total running time of the algorithm is $\Theta(n^4)$.

Another interesting special case occurs when each $\pi_k$ is the cyclic perm
$[k,\,k-1,\,\ldots\,,\,1]$. Here we find that $\sigma_{kj}$ takes
$$ x\mapsto\cases{x-(k-j),&if $x>k-j$;\cr
          k+1-x,&if $x\le k-j$.\cr}$$
It turns out that we have
$$\sigma_{kj}\pi_i=\cases{
\sigma_{k-j,1}\,\sigma_{k-j+1,1}\,\sigma_{k-j+i,k-j+i-1}\,\sigma_{kj}\,,
 &if $i<j$;\cr
\sigma_{k-i,1}\,\sigma_{k-j,1}\,\sigma_{k-j+1,k-i+1}\,\sigma_{k,j-1}\,,
 &if $1<j<i<k$;\cr
\sigma_{k-i,1}\,\sigma_{k-2,i-2}\,\sigma_{k-1,k-i+1}\,\sigma_{ki}\,,
 &if $j=1$ and $2<i<k$;\cr
\sigma_{k2}\,,&if $j=1$ and $i=2$;\cr
\sigma_{k-1,1}\,,&if $j=1$ and $i=k$;\cr
\sigma_{k-j,1}\,\sigma_{k-j+1,1}\,\sigma_{k,j-1}\,,&if $1<j<i=k$.\cr}$$
So the memberships tests need at most 4 multiplications each, and again the total
running time is $\Theta(n^4)$.

In both of these special cases, it turns out that the iterative
implementation of step~A2 will also define the same
perms~$\sigma_{kj}$. Hence the running time will be $\Theta(n^4)$
under either of the implementations we have discussed.

It is interesting to analyze the algorithm in another special case,
when there are just two generators $\sigma_n=[1,2,\ldots,n]$ and
$\tau_n=[n-1,n]$. Assume that the recursive implementation is used.
First, Algorithm $A_n(\sigma_n)$ sets $T(n)=\{\sigma_n\}$ and performs
$B_n(\sigma_n)$. Algorithm $B_n(\sigma_n)$ sets
$\sigma_{n1}\gets\sigma_n$ and performs $B_n(\sigma_n^2)$, which
sets $\sigma_{n2}\gets\sigma_n^2$ and performs
$B_n(\sigma_n^3)$, etc. Thus $\sigma_{nj}$ becomes $\sigma_n^j$ for
all~$j$. Second, Algorithm $A_n(\tau_n)$ adds~$\tau_n$ to $T(n)$ and
performs
$B_n(\tau_n),B_n(\sigma_n\tau_n),\ldots,B_n(\sigma_n^{n-1}\tau_n)$.
The first of these subroutines, $B_n(\tau_n)$, performs algorithm
$A_{n-1}(\tau_n\sigma_n)$, which is $A_{n-1}(\sigma_{n-1})$. The
second subroutine, $B_n(\sigma_n\tau_n)$, performs
$A_{n-1}(\sigma_n\tau_n\sigma_n^{-1})$, which is
$A_{n-1}(\tau_{n-1})$. Therefore we can use induction on~$n$ to show
that $\sigma_{kj}=\sigma_k^j$ for all~$j$ and~$k$. It is easy to
verify that each membership test requires at most three nontrivial
multiplications. Therefore the total running time in this special case
comes to only $\Theta(n^3)$, although $\Gamit(n)$ is the full
symmetric group $\Piit(n)$.

\medskip
\noindent
{\bf 7.\enspace A random example.}\enspace
 The conditions of the construction in Section~6 allow
$(k-1)!$ possibilities for each perm~$\pi_k$. Let us consider the average
total running time of the algorithm when each of the $1!\,2!\ldots(n-1)!$
choices of $\{\pi_2,\pi_3,\ldots,\pi_n\}$ is equally likely.
On intuitive grounds it appears plausible that the average running time
will be $\Theta(n^5)$, because most of the multiplications in a ``random''
situation will be by non-identity perms. This indeed turns out to be true,
at least when the recursive implementation of step~A2 is used;
but the proof is a bit delicate.

As before, the running time is dominated by $\Theta(n^3)$ successful tests for
membership of $\sigma_{kj}\pi_i$ in $\Gamit(k)$, where $k>j\ge1$ and
$k\ge i>1$ and $i\ne j$. We know that the total running time is $O(n^5)$, so
we need only show that the average value is $\Omega(n^5)$; and for this
purpose it will suffice to consider only the membership tests with $k>j>i$.

The membership test for $\sigma_{kj}\pi_i$ performs the multiplications
$$\sigma_{kj}\pi_i\sigma_{kj_k}^-\sigma_{k-1,j_{k-1}}^-\ldots
\sigma_{2j_2}^-\,,$$
and the cost is $l$ for each multiplication such that $j_l\ne l$. Since $j>i$,
we always have $j_k=j$. Let us fix the values $k$, $j$, $i$, and~$l$, where
$k>j>i>1$ and $k>l>i$, and try
to determine an upper bound for the probability that $j_l=l$. The
following analysis applies to any given (not necessarily random)
sequence of perms $\pi_l,\ldots,\pi_2$, with $\pi_k,\ldots,\pi_{l+1}$
varying randomly.

Let $i-r$ be the number of points $\le i$ that are
fixed by the given perm~$\pi_i$. By assumption, $\pi_i$ takes $i\mapsto i-1$,
hence $r\ge2$.

Our first goal is to determine the probability that we have
$j_{k-1}=k-1$, $j_{k-2}=k-2$, \dots,~$j_l=l$. This holds iff
$\sigma_{kj}\pi_i\sigma_{kj}^-\in\Piit(l-1)$. Note that,
in the recursive implementation of step~A2, we have
$$\sigma_{kj}\pi_i\sigma_{kj}^-=
\pi_k\pi_{k-1}\ldots\pi_{j+1}\pi_i\pi_{j+1}^-\ldots\pi_{k-1}^-\pi_k^-
=\pi_k\rho\pi_k^-,$$
where $\rho$ is a perm of $\Piit(k-2)$ that has the same cycle structure
as $\pi_i$; hence $\rho$ fixes exactly $k-2-r$ points $\le k-2$. Consider
what happens to $\pi_k\rho\pi_k^-$ as $\pi_k$ runs through its $(k-1)!$
possible values: We obtain a uniform distribution
over all perms of $\Piit(k-1)$ having the same cycle structure as $\rho$.
For example, if $r=7$ and $\rho=[1\,2\,7][3\,6][4\,9]$, the $(k-1)!$
perms $\pi_k\rho\pi_k^-$ are just $[a_1\,a_2\,a_7][a_3\,a_6][a_4\,a_9]$ as
$a_1\ldots a_{k-1}$ runs through the images of all perms of $\Piit(k-1)$.
Therefore the probability that $\sigma_{kj}\pi_i\sigma_{kj}^-\in\Piit(l-1)$
is
$${l-1\choose r}\bigg/{k-1\choose r}={(l-1)(l-2)\ldots(l-r)\over
 (k-1)(k-2)\ldots(k-r)}\,.$$

Now let's compute the probability that $j_{k-1}=k-1$, \dots, $j_{q+1}=q+1$,
$j_q<q$, and $j_l=l$, given a subscript~$q$ in the range $k>q>l$. We will
assume that $\pi_{k-1}$, \dots, $\pi_{q+1}$, $\pi_{q-1}$, \dots,~$\pi_2$
have been assigned some fixed values, while $\pi_k$ and $\pi_q$ run
independently through all of their $(k-1)!\,(q-1)!$ possibilities. Under
these circumstances we will prove that $\sigma_{kj}\pi_i\sigma_{kj}^-
\sigma_{qj_q}^-$ is uniformly distributed over $\Piit(q-1)$.

Let $p$ be a positive integer less than $q$. Let $\alpha\in\Piit(q)$
take $q\mapsto p$ and have the same cycle structure as $\pi_i$. Also let
$\beta$ be an element of $\Piit(q-1)$. Then there is exactly one perm
$\pi_q$ that will make $\alpha\sigma_{qp}^-=\beta$, namely
$$\pi_q=\beta^-\alpha\pi_{p+1}^-\ldots\pi_{q-1}^-\,.$$
(This perm takes $q\mapsto q-1$ and fixes all points $>q$, so it meets
the conditions necessary to be called $\pi_q$.) Moreover, when $\pi_q$
has this value, the number of perms $\pi_k$ such that $\sigma_{kj}\pi_i
\sigma_{kj}^-=\alpha$ is independent of $\alpha$, as we have observed in
the previous case. Therefore the probability that ($j_q=p$ and $\sigma_{kj}
\pi_i\sigma_{kj}^-\sigma_{qp}^-=\beta$) is independent of~$\beta$,
and independent of~$p$.

The uniform distribution of $\sigma_{kj}\pi_i\sigma_{kj}^-\sigma_{qp}^-$
implies that we have ($j_{k-1}=k-1$, \dots, $j_{q+1}=q+1$, $j_q<q$, and
$j_l=l$) with probability $1/l$ times the probability that ($j_{k-1}=k-1$,
\dots, $j_{q+1}=q+1$, and $j_q<q$), because the values $j_{q-1}\ldots j_2$
are uniformly distributed. And we know from the previous analysis that this
probability is
$${1\over l}\left(q(q-1)\ldots(q-r+1)-(q-1)(q-2)\ldots(q-r)\over
 (k-1)(k-2)\ldots(k-r)\right)={r\over l}\,{(q-1)\ldots(q-r+1)\over
 (k-1)\ldots(k-r)}\,.$$

Finally, therefore, we can compute the probability that $j_l=l$, when
$k$, $j$, $i$, and $l$ are given as above and $\pi_i$ has $i-r$ fixed
points: It comes to
$$\openup1\jot
\eqalign{&{1\over(k-1)\ldots(k-r)}\bigg((l-1)\ldots(l-r)+{r\over l}\sum_{l<q<k}
 (q-1)\ldots(q-r+1)\biggr)\cr
&\qquad={1\over l}+{(l-1)\ldots(l-r)(l-r-1)\over(k-1)\ldots(k-r)(l-r)}\cr
&\qquad<{1\over l}+{(l-1)\ldots(l-r)\over(k-1)\ldots(k-r)}\,.\cr}$$
Since $r\ge2$, we obtain the desired upper bound
$$\Pr(j_l=l)<{1\over l}+{(l-1)(l-2)\over(k-1)(k-2)}<{1\over l}+{l^2\over k^2}\,.$$
This implies the desired lower bound $\Omega(n^5)$ on the total multiplication time.
We can, for example, sum over $\Omega(n^4)$ values $(k,j,i,l)$ with
$1<i\le{1\over4}n<l\le{1\over2}n<j\le{3\over4}n<k\le n$; in each of these
cases a multiplication will require $\Omega(n)$ steps with probability
at least $1-(1/l+l^2/k^2)>1/2$ when $n\geq 72$.

Since the average running time is $\Omega(n^5)$, there must exist,
for all~$n$, a sequence of perms $\pi_2$,~\dots,~$\pi_n$ that make the
algorithm do $\Omega(n^5)$ operations. But it appears to be difficult to
define such perms via an explicit construction.
Nor is there an obvious way to prove the $\Omega(n^5)$ bound when the
iterative implementation of step~A2 is adopted in place of the
recursive implementation, even in the totally random case.

\medskip
\noindent
{\bf 8.\enspace More meaningful upper bounds.}\enspace
The examples studied above show that it
is misleading to characterize algorithms $A$ and~$B$ by merely saying that
they will process $m$~perms of $\Piit(n)$ with a
 worst-case running time of $O(n^5+mn^2)$.
In one sense this estimate is sharp, because we've seen that
$\Omega(n^5)$ behavior may indeed occur; but our other examples, together
with extensive computational experience, show that the procedures often
run considerably faster in practice.

We can improve the estimate of Section 4 by introducing another
parameter. Let $g$ be the order of the group $\Gamit(n)$ that is generated.
Then we have the following result:

\proclaim Theorem. A transversal system for a perm group of order~$g$
generated by $m$~perms of~$\Piit(n)$ can be found in at most\/
$O\bigl(n^2(\log g)^3/\!\log n\bigr)+O\bigl(n^2(\log g)^2\bigr)
+O(mn\log g)$ steps, using at most $O(n^2\log g/\!\log
n)+O\bigl(n(\log g)^2\bigr)$ memory cells.

\noindent Proof. Let $s(k)$ and $t(k)$ be defined as before. Then
$g=\prod_{k=1}^n s(k)$, and the number of membership tests is $m+
\sum_{k=1}^n\bigl(s(k)t(k)-s(k)+1\bigr)$. Each membership test involves
at most $O(\log g)$ multiplications by non-identity perms, because the
number of indices~$k$ with $s(k)>1$ cannot exceed $\theta(g)$, the
total number of prime factors of~$g$. This accounts for the term
$O(mn\log g)$ in the theorem.

Moreover, each $t(k)$ is at most $\theta(g)=O(\log g)$, as we have
argued before. Therefore we can complete the proof of the time bound
by showing that
$\sum_{k=1}^n\bigl(s(k)-1\bigr)=O(n\log g/\!\log n)$.

Given $n$ and $s$, let us try to minimize the product $\prod_{k=1}^n s_k$
subject to the conditions
$$s=\sum_{k=1}^n(s_k-1)\qquad\hbox{and}\qquad 1\le s_k\le k\,.$$
If $s_{k-1}>s_k$, we can interchange $s_{k-1}\leftrightarrow s_k$
without violating the conditions; hence we may assume that $s_1\le s_2
\le\cdots\le s_n$. Furthermore, if $1<s_{k-1}\le s_k<k$, we can decrease
the product by setting $(s_{k-1},s_k)\gets(s_{k-1}-1,s_k+1)$. Hence the
product is smallest when we have $s_k=k$ for as many large~$k$ as possible:
$$s_n=n,\quad s_{n-1}=n-1,\quad \ldots,\quad s_{q+1}=q+1,\quad s_q=r,\quad
 s_{q-1}=\cdots=s_1=1.$$
Here $q$ and $r$ are the unique integers such that
$$\textstyle{n\choose2}-s-1\,=\,{q\choose2}-r\qquad\hbox{and}\qquad 1\le r<q\le n.$$
(We assume that $0\le s<{n\choose2}$.) The minimum product is
$$P(n,s)=r{n!\over q!}\,.$$
The actual product in the algorithm is $g\ge P\bigl(n,\sum_{k=1}^n\bigl(
s(k)-1\bigr)\bigr)$, hence our proof will be complete if we can show that
$$s=O\Bigl(n{\log P(n,s)\over\log n}\Bigr)\,.$$
But this is not difficult. If $s\ge{1\over4}n^2$ we have $q\le n/\sqrt2$,
hence $\log P(n,s)=\Theta(n\log n)$ and the result holds. At the other
extreme, if $0\le s<n$, we have $P(n,s)=s+1$ and again the result is trivial.
Otherwise we note that $n-q\ge\lfloor s/n\rfloor$, hence
$$P(n,s)\ge{n!\over q!}>q^{\lfloor s/n\rfloor}>\Bigl({n\over2}\Bigr)^{\!
s/n-1};$$
the relation $(s/n)\log n=O\bigl(\log P(n,s)\bigr)$ follows immediately.

The space required to store the transversal perms $\sigma_{kj}$ is
$\sum_{k=1}^nk\bigl(s(k)-1\bigr)=O(n^2\log g/\!\log n)$. The space
required to store the strong generators can be reduced to $\sum
k\,t(k)$ summed over those~$k$ with $s(k)>1$, for if $s(k)=1$ we have
$T(k)=T(k-1)$. This sum has $O(\log g)$ terms, each of which is
$O(n\log g)$. So the proof of the theorem is complete.

Inspection of this proof shows that the running time is actually bounded
by a slightly smaller estimate than claimed, namely
$$O\bigl(n^2 l_n(g)^2 \log_n g\bigr)+
O\bigl(n^2l_n(g)^2\bigr)+O\bigl(mnl_n(g)\bigr)\,,
\qquad\hbox{where $l_n(g)=\min\bigl(n,\theta(g)\bigr)$}.$$
The space bound is, similarly, $O(n^2\log_ng)+O\bigl(nl_n(g)^2\bigr)$.
And the examples in Sections 5 and~6 above show that even this improved
bound might be unduly pessimistic; sometimes a judicious relabeling of
points will speed things up.

The storage occupied by strong generators is usually less than the
storage required for perms of the traversal system, but it can be
greater. For example, when $n$ is even and the generators are respectively
$$\vcenter{\halign{$\hfil#\hfil$\cr
[n-1,n]\cr
[n-3,n-2]\,[n-1,n]\cr
\vdots\cr
[1,2]\,\ldots\,[n-3,n-2]\,[n-1,n]\cr}}$$
then $g=2^{n/2}$ and the $n\,l_n(g)^2$ term dominates.

The values of $l_n(g)$ and $\log_n g$ are often substantially smaller than~$n$,
in perm groups of computational interest. For example, the Hall-Janko group
$J_2$ has $g=2^7\cdot3^3\cdot5^2\cdot7$ and $n=100$ (see~[6]); here
$\theta(g)=13$ and $\log_n g\approx2.9$. The unitary group
$U_6(2)$, which has order $g=2^{15}\cdot3^6\cdot5\cdot7\cdot11$, is
represented as a perm group on $n=672$ points in 
the Cayley library (see~[10]);
 in this case $l_n(g)=24$ and $\log_n g\approx3.5$.
Some representative large examples are Conway's perfect group~$\cdot0$,
for which $g=2^{22}\cdot3^9\cdot5^4\cdot7^2\cdot11\cdot13\cdot23$, $n=
196560$, and $\log_n g\approx3.6$; and Fischer's simple group
$F'_{24}$, for which $g=2^{21}\cdot3^{16}\cdot5^2\cdot7^3\cdot11\cdot13\cdot
17\cdot23\cdot29$, $n=306936$, and $\log_n g\approx4.4$. (See~[3].)

\medskip
\noindent
{\bf 9.\enspace Historical remarks and acknowledgments.}\enspace The algorithm
described above is a variant of a fundamental procedure sketched by
Sims in 1967 [8], which he described more fully a few years later 
as part of a larger body of algorithms~[9].
The principal difference between the method of~[9] and the present
method is that Sims essentially worked with sets of strong generators
satisfying the condition $T(1)\subseteq
T(2)\subseteq\,\cdots\,\subseteq
T(n)$. Thus, for example, when $\sigma\in\Sigit(n)$ he would test
the product $\sigma\tau$ for all strong generators~$\tau$; the present
algorithm tests $\sigma\tau$ for such~$\sigma$ only with the
perms~$\tau$ of $T(n)$, namely the given generators~$\pi$.
His example, in which the group generated by $[1,2,4,5,7,3,6]$ and
$[2,4]\,[3,5]$ required the verification of 54 products~$\sigma\tau$,
requires the testing  of only 40~products in the present scheme.
On the other hand, his method for representing the $\Sigit(k)$ as
words in the generators was considerably more economical in its use of
storage space, and space was an extremely critical resource at the
time. Moreover, his way of maintaining strong generators blended well
with the other routines in his system, so it is not clear that he
would have regarded the methods of the present paper as an
improvement.

Polynomial bounds on the worst-case running time were not obvious
from this original work.
Furst, Hopcroft, and Luks showed in 1980 [5] that a transversal system
and a set of strong generators could be found in $O(n^6)$ steps. (In
their method the transversal system and strong generators were
identical.) The author developed the present algorithm a year later,
while preparing to write Volume~4 of
{\sl The Art of Computer Programming} and while advising 
Eric~W. Hamilton, an undergraduate
student who was working on a research project with Persi Diaconis~[4].
 The present method became more widely known after
the author discussed it informally at a conference in Oberwolfach on
November~6, 1981; several people, notably Clement Lam, suggested
clarifications of the rough notes that were distributed at that time.
Eventually Professor Babai was kind enough to suggest that the notes of
1981 be published now, instead of waiting until Volume~4 has been completed.
Those notes are reproduced with slight improvements in Sections 1--4 of the
present paper. The author is grateful to the referees and to Profs.\ Babai and Luks
for several penetrating remarks that prompted the additional material in
Sections 5--8.

Improved methods have been discovered in the meantime, notably by Jerrum~[7],
who has reduced the worst-case storage requirement to order~$n^2$. 
Babai, Luks, and Seress~[2]
have developed a more complicated procedure whose worst case running time is
only $O(n^{4+\epsilon})$.

The word ``perm,'' introduced experimentally in the author's Oberwolfach notes,
does not seem to be winning any converts. (In fact, Pratt himself has forgotten
that he once made this suggestion in conversation with the author.) However,
the proposal to use the notation $\pi^-$ for inverses, instead of the
usual $\pi^{-1}$, has significantly greater merit, and the author hopes to
see it widely adopted in future years.
The shorter notation is easier to write on a blackboard and easier to type
on a keyboard. Moreover, the longer
notation $\alpha^{-1}$ is redundant, just as $\alpha^1$ is
redundant; in fact, $\alpha^{-1}$ stands for $\alpha^-$ raised to the first
power! Thus there is no conflict between the two conventions, and a
gradual changeover should be possible.

\bigskip
\centerline{\bf References}

\medskip
\disleft 20pt:[1]:
L\'aszl\'o Babai,
``On the length of subgroup chains in the symmetric group,''
{\sl Communications in Algebra\/ \bf 14}
(1986), 1729--1736.

\medskip
\disleft 20pt:[2]:
L\'aszl\'o Babai, Eugene M. Luks, and \'Akos Seress,
``Fast management of permutation groups,''
{\sl 29th Annual Symposium on Foundations of Computer Science\/}
(IEEE Computer Society, 1988), 272--282.

\medskip
\disleft 20pt:[3]:
J. H. Conway,
``Three lectures on exceptional groups,''
in M. B. Powell and G. Higman, ed.,
{\sl Finite Simple Groups}, Proceedings of the Oxford Instructional Conference
  on Finite Simple Groups, 1969
(London: Academic Press, 1971), 215--247.

\medskip
\disleft 20pt:[4]:
Persi Diaconis, R. L. Graham, and William M. Kantor,
``The mathematics of perfect shuffles,''
{\sl Advances in Applied Mathematics\/ \bf 4}
(1983), 175--196.

\medskip
\disleft 20pt:[5]:
Merrick Furst, John Hopcroft, and Eugene Luks,
``Polynomial-time algorithms for permutation groups,''
{\sl 21st Annual Symposium on Foundations of Computer Science\/}
(IEEE Computer Society, 1980), 36--41.

\medskip
\disleft 20pt:[6]:
Marshall Hall, Jr.\ and David Wales,
``The simple group of order 604,800,''
{\sl Journal of Algebra\/ \bf 9}
(1968), 417--450.

\medskip
\disleft 20pt:[7]:
Mark Jerrum,
``A compact representation for permutation groups,''
{\sl Journal of Algorithms\/ \bf7}
(1986), 60--78.

\medskip
\disleft 20pt:[8]:
Charles C. Sims,
``Computational methods in the study of permutation groups,''
in John Leech, ed., {\sl Computational Problems in Abstract Algebra},
Proceedings of a conference held at Oxford University in 1967
(Oxford: Pergamon, 1970), 169--183.

\medskip
\disleft 20pt:[9]:
Charles C. Sims,
``Computation with permutation groups,'' in S.~R. Petrick, ed., 
{\sl Proc.\ Second Symposium on Symbolic and Algebraic Manipulation},
Los Angeles, California (New York: ACM, 1971), 23--28.

\medskip
\disleft 20pt:[10]:
D. E. Taylor,
``Pairs of generators for matrix groups,''
{\sl The Cayley Bulletin\/ \bf 3} (Department of Pure Mathematics,
University of Sydney, 1987).

\bye